\documentclass[a4paper,reqno]{amsart}

\usepackage{amssymb,amsmath}

\usepackage{paralist,enumitem}
\setlist{listparindent=0pt,parsep=3pt}
\setdefaultenum{1.}{}{}{}

\usepackage[pdftex,linktocpage,pdfstartview=FitH,colorlinks=true,allcolors=blue]{hyperref}
\usepackage{amsthm} 
\usepackage[capitalise]{cleveref}

\usepackage[initials,msc-links]{amsrefs}[2007/10/22]

\newtheorem{thm}{Theorem}[section]

\newtheorem*{proposition*}{Proposition}
\newtheorem*{theorem*}{Theorem}
\theoremstyle{definition}

\theoremstyle{remark}
\newtheorem*{rem}{Remark}

\numberwithin{equation}{section}

\newcommand{\R}{\mathbb{R}}

\renewcommand{\P}{\mathbb{P}}
\renewcommand{\H}{\mathbb{H}}

\DeclareMathOperator{\im}{im}

\DeclareMathOperator{\End}{End}

\newcommand{\Span}[1]{\langle#1\rangle}

\DeclareSymbolFont{script}{U}{eus}{m}{n}
\DeclareSymbolFontAlphabet{\mathscr}{script}
\DeclareMathSymbol{\EuWedge}{0}{script}{"5E}
\newcommand{\Wedge}{\EuWedge}
\renewcommand{\d}{\mathop{}\!\mathrm{d}}

\DeclareMathOperator{\rO}{O}

\newcommand{\fso}{\mathfrak{so}}

\renewcommand{\q}{\mathfrak{q}}

\newcommand{\Ad}{\operatorname{Ad}}

\newcommand{\II}{\mathrm{I\kern-0.5ptI}}

\newcommand{\sect}[1]{\Gamma\,#1}

\newcommand{\sub}{\subseteq}
\newcommand{\st}{\mathrel{|}}
\newcommand{\set}[1]{\{#1\}}
\newcommand{\vol}{\mathrm{vol}}

\renewcommand{\Re}{\operatorname{Re}}
\renewcommand{\Im}{\operatorname{Im}}

\renewcommand{\lor}[1][4,1]{\R^{#1}}
\newcommand{\cL}{\mathcal{L}}
\newcommand{\ip}[1]{(#1)}

\newcommand{\half}{\tfrac12}

\newcommand{\TitleWithUrl}[1]{\IfEmptyBibField{doi}%
  {\IfEmptyBibField{url}{\textit{#1}}%
    {\IfEmptyBibField{eprint}{\href {\BibField{url}}{\textit{#1}}}{\textit{#1}}}%
    }%
  {\href {https://doi.org/\BibField{doi}}{\textit{#1}}}}
\renewcommand{\eprint}[1]{\IfEmptyBibField{url}{\url{#1}}%
  {\href {\BibField{url}}{#1}}}

\BibSpec{article}{%
    +{}  {\PrintAuthors}                {author}
    +{,} { \TitleWithUrl}               {title}
    +{.} { }                            {part}
    +{:} { \textit}                     {subtitle}
    +{,} { \PrintContributions}         {contribution}
    +{.} { \PrintPartials}              {partial}
    +{,} { }                            {journal}
    +{}  { \textbf}                     {volume}
    +{}  { \PrintDatePV}                {date}
    +{,} { \issuetext}                  {number}
    +{,} { \eprintpages}                {pages}
    +{,} { }                            {status}
    +{,} { available at \eprint}        {eprint}
    +{}  { \parenthesize}               {language}
    +{}  { \PrintTranslation}           {translation}
    +{;} { \PrintReprint}               {reprint}
    +{.} { }                            {note}
    +{.} {}                             {transition}
    +{}  {\SentenceSpace \PrintReviews} {review}
}

\BibSpec{collection.article}{%
    +{}  {\PrintAuthors}                {author}
    +{,} { \TitleWithUrl}                     {title}
    +{.} { }                            {part}
    +{:} { \textit}                     {subtitle}
    +{,} { \PrintContributions}         {contribution}
    +{,} { \PrintConference}            {conference}
    +{}  {\PrintBook}                   {book}
    +{,} { }                            {booktitle}
    +{,} { \PrintDateB}                 {date}
    +{,} { pp.~}                        {pages}
    +{,} { }                            {status}
    +{,} { available at \eprint}        {eprint}
    +{}  { \parenthesize}               {language}
    +{}  { \PrintTranslation}           {translation}
    +{;} { \PrintReprint}               {reprint}
    +{.} { }                            {note}
    +{.} {}                             {transition}
    +{}  {\SentenceSpace \PrintReviews} {review}
}
\BibSpec{book}{%
    +{}  {\PrintPrimary}                {transition}
    +{,} { \TitleWithUrl}               {title}
    +{.} { }                            {part}
    +{:} { \textit}                     {subtitle}
    +{,} { \PrintEdition}               {edition}
    +{}  { \PrintEditorsB}              {editor}
    +{,} { \PrintTranslatorsC}          {translator}
    +{,} { \PrintContributions}         {contribution}
    +{,} { }                            {series}
    +{,} { \voltext}                    {volume}
    +{,} { }                            {publisher}
    +{,} { }                            {organization}
    +{,} { }                            {address}
    +{,} { \PrintDateB}                 {date}
    +{,} { }                            {status}
    +{}  { \parenthesize}               {language}
    +{}  { \PrintTranslation}           {translation}
    +{;} { \PrintReprint}               {reprint}
    +{.} { }                            {note}
    +{.} {}                             {transition}
    +{}  {\SentenceSpace \PrintReviews} {review}
}

\BibSpec{thesis}{%
    +{}  {\PrintAuthors}                {author}
    +{.} { \PrintDate}                  {date}
    +{.} { \TitleWithUrl}               {title}
    +{:} { \textit}                     {subtitle}
    +{,} { \PrintThesisType}            {type}
    +{,} { }                            {organization}
    +{,} { }                            {address}
    +{,} { \eprint}                     {eprint}
    +{,} { }                            {status}
    +{}  { \parenthesize}               {language}
    +{}  { \PrintTranslation}           {translation}
    +{;} { \PrintReprint}               {reprint}
    +{.} { }                            {note}
    +{.} {}                             {transition}
    +{}  {\SentenceSpace \PrintReviews} {review}
}

\title{Bonnet pairs, isothermic surfaces and the retraction form}

\author{F.E. Burstall}
\address{Department of Mathematical Sciences\\ University of Bath\\
  Bath BA2 7AY\\UK}
\email{feb@maths.bath.ac.uk}

\author{T. Hoffmann}
\address{Department of Mathematics\\ Technical University of Munich\\
  85748 Garching\\Germany}
\email{tim.hoffmann@ma.tum.de}

\author{F. Pedit}
\address{Department of Mathematics and Statistics\\ University of Massachusetts Amherst\\
  Amherst, MA 01030\\USA}
\email{pedit@math.umass.edu}

\author{A.O. Sageman-Furnas}
\address{Department of Mathematics\\ North Carolina State University\\
  Raleigh, NC 27607\\USA}
\email{asagema@ncsu.edu}

\begin{document}

\begin{abstract}
  We give a modern account of the classical theory of
  Bianchi \cite{Bia03} (see also \cite{KamPedPin98})
  relating isothermic surfaces to Bonnet pairs.  The main
  novelty is to identify the derivatives of the Bonnet pair
  with a component of the retraction form of the isothermic
  surface.
\end{abstract}
\maketitle

\section{Isothermic surfaces}
\label{sec:isothermic-surfaces}

We begin by briefly rehearsing the notion of an isothermic
surface in the conformal $3$-sphere and its retraction form:
a closed $1$-form on $\Sigma$ with values in the Lie algebra
of conformal vector fields on $S^3$.

\subsection{Light-cone model}
\label{sec:light-cone-model}

View the conformal $3$-sphere as the projective light-cone
$\P(\cL)$.  Thus $\lor$ is a $5$-dimensional vector space with
inner product $\ip{\,\,,\,}$ of signature $(4,1)$ and
$\cL=\set{v\in\lor\st\ip{v,v}=0}$.  The natural action of
$\rO(4,1)$ on $\P(\cL)$ is by conformal diffeomorphisms so
that $\rO(4,1)$ is a double cover of the conformal
diffeomorphism group of $S^3$.

In particular, we identify $\fso(4,1)$ with the Lie algebra
of conformal vector fields of $S^3$.  We further identify
$\fso(4,1)$ with $\Wedge^2\lor$ via
\begin{equation}
  \label{eq:1}
  (a\wedge b)c=\ip{a,c}b-\ip{b,c}a.
\end{equation}

Let $f\colon\Sigma\to S^3=\P(\cL)$ be an immersion of an
oriented surface into the conformal $3$-sphere.  Recall that
$f$ is \emph{(globally) isotermic} if there is a non-zero
holomorphic quadratic differential $q$ (for the Riemann
surface structure on $\Sigma$ induced by $f$) which commutes
in an appropriate way with the second fundamental form of
$f$.  This structure is encapsulated in an
$\fso(4,1)$-valued $1$-form $\eta$ in the following way.

First, view $f$ as a null line subbundle
$f\leq\underline{\R}^{4,1}:=\Sigma\times\lor$ and then
define the \emph{retraction form}
$\eta\in\Omega^{1}_{\Sigma}(\Wedge^{2}\lor)$ of $f$ by
\begin{equation}
  \label{eq:2}
  \eta=(\d\sigma\circ Q^{\#}_{\sigma})\wedge\sigma
\end{equation}
where $\sigma\in\sect f^{\times}$ (thus
$\sigma\colon\Sigma\to\cL\sub\lor$ lifts $f$) and
$Q^{\#}_{\sigma}\in\sect\End(T\Sigma)$ is given by
\begin{equation*}
  2\Re q=: Q = \ip{\d\sigma\circ Q^{\#}_{\sigma},\d\sigma}.
\end{equation*}
It is easy to check that $\eta$ is independent of choice of
lift $\sigma$.

The essential point is that the retraction form and, indeed,
isothermicity of the surface is
completely characterised by two properties:
\begin{thm}[c.f.\ \cite{BurSan12}*{Proposition~1.4~}]
  \label{th:1}
  A surface $f$ is isothermic if and only if it admits a
  non-zero $1$-form $\eta\in\Omega^1_{\Sigma}(\Wedge^{2}\lor)$ such
  that
  \begin{compactenum}
  \item $\eta$ takes values in $f\wedge f^{\perp}$;
  \item $\d\eta=0$.
  \end{compactenum}
  In this case, the holomorphic quadratic differential is
  recovered by
  \begin{equation*}
    Q\sigma=\eta\d\sigma,
  \end{equation*}
  for $\sigma\in\sect f^{\times}$ and \eqref{eq:2} holds.
\end{thm}

The essential point about this construction is that $\eta$
is manifestly M\"obius-invariant: for $g\in\rO(4,1)$, $gf$
is isothermic with retraction form $\Ad_g\eta$.
\subsection{Spherical model}
\label{sec:spherical-model}

Now fix $\q\in\lor$ with $\ip{\q,\q}=-1$ and decompose
$\lor$ into an orthogonal direct sum
\begin{equation*}
  \lor=\R^4\oplus\Span{\q}.
\end{equation*}
We get a corresponding decomposition
\begin{equation}\label{eq:3}
  \Wedge^2\lor=\Wedge^2\R^4\oplus\R^4\wedge\Span{\q}
\end{equation}
with $\Wedge^2\R^4\cong\fso(4)$ via \eqref{eq:1}.

Let $f\colon\Sigma\to\P(\cL)$ be an isothermic surface with
retraction form $\eta$ and let $y\in\sect f^{\times}$ be the
unique lift with $\ip{y,\q}=-1$.  Then $y=x+\q$ for
$x\colon\Sigma\to S^3\sub\R^4$ and \eqref{eq:2} reads:
\begin{equation}
  \label{eq:4}
  \eta=\omega\wedge(x+\q)=\omega\wedge x +\omega\wedge\q
\end{equation}
where $\omega=\d x\circ
Q^{\#}_{y}\in\Omega_{\Sigma}^1(\R^4)$.  Clearly, $\d\eta=0$
if and only if both $\omega$ and $\omega\wedge x$ are
closed.  However,
\begin{equation*}
  \d(\omega\wedge x)=\d\omega\wedge x-\omega\curlywedge\d x,
\end{equation*}
where $\curlywedge$ is exterior product of $\R^4$-valued
$1$-forms using wedge product in $\R^4$ to multiply
coefficients.  We conclude that
the closure of $\eta$ amounts to
\begin{subequations}\label{eq:7}
  \begin{align}
    \label{eq:5}
    \d\omega&=0\in\Omega^2_{\Sigma}(\R^4)\\
    \label{eq:6}
    \omega\curlywedge\d x&=0\in\Omega^2_{\Sigma}(\Wedge^2\R^4).
  \end{align}
\end{subequations}

\begin{rem}
  One way to understand these equations is to view $x$ as a
  surface in $\R^4$.  Then $x$ is an isothermic surface in
  $\R^4$ with Christoffel dual\footnote{Thus $x$ and $x^{*}$
    have parallel tangent planes, the same conformal
    structure induced on $\Sigma$ with opposite
    orientations.} $x^{*}$ given by (locally) integrating
  $\omega=\d x^{*}$.  Remark that $x^{*}$ immerses off the
  zero-set of $\eta$ which is the zero divisor of $q$.
\end{rem}

\section{Bonnet pairs}
\label{sec:bonnet-pairs}

Surfaces $F_{\pm}\colon\Sigma\to\R^3$ are said to be a
\emph{Bonnet pair} if they are isometric and, for a suitable
choice of unit normals $n_{\pm}$, their mean curvatures
coincide.  Equivalently, their second fundamental forms
$\II_{\pm}$ differ by a trace-free symmetric bilinear form
which, thanks to the Codazzi equation, is of the form
$\Re\hat{q}$, for some holomorphic quadratic differential
$\hat{q}$.

With $x\colon\Sigma\to S^3\sub\R^4$ as above, suppose that
$\omega\wedge x$ is exact and $\eta$ is never zero so that
we have an immersion
$F\colon\Sigma\to\fso(4)$ with $\d F=\omega\wedge x$. There
is a well-known Lie algebra decomposition
\begin{equation}
  \label{eq:8}
  \fso(4)=\fso(3)\oplus\fso(3)
\end{equation}
and so we write $F=F_++F_-$ with
$F_{\pm}\colon\Sigma\to\fso(3)\cong\R^{3}$.

With this notation established, here is our formulation of
Bianchi's result:

\begin{thm}\label{th:2}
  Let $f=\Span{x+\q}$ be isothermic with never-zero
  holomorphic quadratic differential $q$ and
  $\eta=\omega\wedge\q+\omega\wedge x$.  Suppose that
  $\omega\wedge x=\d F$ and write $F=F_++F_-$ as above.
  Then $F_{\pm}$ are a Bonnet pair:
  \begin{compactenum}
  \item $F_{\pm}$ are isometric and are conformal to $f$.
  \item $\II_+-\II_{-}=2\sqrt{2}\Re(iq)$.
  \end{compactenum}

  Moreover, up to a sign, all Bonnet pairs arise this way:
  if $F_{\pm}$ are a Bonnet pair with $\II_{+}-\II_-$ never
  zero, then one of $F_+\pm F_{-}\colon\Sigma\to\fso(4)$ has
  derivative $\omega\wedge x$ for
  $x\colon\Sigma\to S^3\sub\R^4$ isothermic.
\end{thm}

For the proof, we start 
by fixing a unit length
$\vol\in\Wedge^4\R^4$ to orient $\R^{4}$ and introduce the Hodge star operator
$S\in\End(\Wedge^2\R^4)$ by
\begin{equation*}
  \ip{\alpha,\beta}\vol=\alpha\wedge S(\beta).
\end{equation*}
Then $S$ is an involutive isometry whose $\pm1$-eigenspaces
are the self-dual and anti-self-dual $2$-vectors
$\Wedge^2_{\pm}\R^{4}$.  Under the identification
\eqref{eq:1}, the decomposition \eqref{eq:8} becomes the (orthogonal)
eigenspace decomposition
\begin{equation*}
  \Wedge^2\R^4=\Wedge^2_+\R^4\oplus\Wedge^2_-\R^4.
\end{equation*}
Let $\kappa$ be the signature $(3,3)$ \emph{Klein inner product} on
$\Wedge^2\R^4$ given by
\begin{equation*}
  \kappa(\alpha,\beta)=\ip{\alpha,S\beta}
\end{equation*}
so that
\begin{equation*}
  \kappa(\alpha,\beta)\vol=\alpha\wedge\beta.
\end{equation*}
Note that a $2$-vector is decomposable if and only if it is
isotropic for $\kappa$.

With all this in hand, let $f$ be isothermic as above with
$\d F=\omega\wedge x$ so that $F_{\pm}=\half(F\pm SF)$.  We
have
\begin{equation*}
  0=\kappa(\d F,\d F)=\ip{\d F,S\d F}=\ip{\d F_+,\d
    F_+}-\ip{\d F_-,\d F_-}
\end{equation*}
so that $F_{\pm}$ are isometric.  Moreover,
\begin{equation*}
  2\ip{\d F\pm,\d F_{\pm}}=\ip{\d F,\d F}=
  \ip{(\d x\circ Q^{\#})\wedge x,(\d x\circ Q^{\#})\wedge x}=
  \ip{(\d x\circ Q^{\#}),(\d x\circ Q^{\#})},
\end{equation*}
since $\ip{x,x}=1$ so that $\ip{\d x,x}=0$.  Since $Q^{\#}$
is symmetric and trace-free, so conformal, we conclude that
the common metric on $F_{\pm}$ is conformal to that of $x$.

Now choose the unit normal $n$ to $x$ in $S^3$ for which
\begin{equation*}
  \vol_x\wedge x\wedge n=\vol,
\end{equation*}
where $\vol_x$ is the volume form of $x$ on the oriented
surface $\Sigma$.  Then $n\wedge x$ is a unit normal to $F$:
\begin{equation}\label{eq:9}
  \ip{\d F,n\wedge x}=\ip{\d x\circ Q^{\#},n}=0.
\end{equation}
Set $n_{\pm}=\tfrac1{\sqrt 2}(n\wedge x\pm S(n\wedge x))$.
Then
\begin{equation*}
  \ip{n_+,n_{+}}+\ip{n_-,n_-}=2
\end{equation*}
while
\begin{equation*}
  0=2\kappa(n\wedge x,n\wedge x)=\ip{n_+,n_{+}}-\ip{n_-,n_-}
\end{equation*}
so that $n_{\pm}$ have unit length.  Moreover,
\eqref{eq:9} reads
\begin{equation*}
  \ip{\d F_{+},n_{+}}+\ip{\d F_-,n_-}=0.
\end{equation*}
On the other hand,
\begin{equation*}
  0=\sqrt{2}\kappa(\d F,n\wedge x)=\ip{\d F_+,n_+}-\ip{\d F_-,n_{-}},
\end{equation*}
and we conclude that $n_{\pm}$ are unit normals to
$F_{\pm}$.

Finally,
\begin{equation*}
  \II_+-\II_-=-\ip{\d F_+,\d n_+}+\ip{\d F_-,\d
    n_-}=-\sqrt{2}\kappa(\d F,\d(n\wedge x))
\end{equation*}
so that
\begin{equation*}
(\II_+-\II_-)\vol=-\sqrt{2}\d F\wedge \d(n\wedge
x)=-\sqrt{2}\omega\wedge x\wedge\d(n\wedge
x)=-\sqrt{2}\omega x\wedge n\wedge\d x.
\end{equation*}
Let $J$ be the orthogonal almost complex structure on
$\Span{x,n}^{\perp}$ for which $U\wedge JU\wedge x\wedge
n>0$ so that $J\d x=\d x J^{\Sigma}$.  A short computation
gives
\begin{equation}
  \label{eq:10}
  \omega\wedge x\wedge n\wedge\d x=-\ip{\omega,J\d x}\vol
\end{equation}
so that
\begin{equation*}
  \II_+-\II_-=\sqrt{2}\ip{\omega,J\d x}=\sqrt{2}\ip{\d
    x\circ Q^{\#},\d x\circ J^{\Sigma}}=2\sqrt{2}\Re(iq).
\end{equation*}
In particular, $F_{\pm}$ are a Bonnet pair.

For the converse, given a Bonnet pair
$F_{\pm}\colon\Sigma\to\Wedge^2_{\pm}\R^{4}$ with unit
normals $n_{\pm}$, define a holomorphic quadratic
differential $q$ by $\II_+-\II_-=2\sqrt{2}\Re(iq)$ which we
assume never vanishes.  Define subbundles $W_{\pm}$ of the
trivial $\Wedge^2\R^{4}$ bundle by
\begin{equation*}
  W_{\pm}=\Span{\im(\d F_+\pm\d F_-),n_+\pm n_-}.
\end{equation*}
Since $F_{\pm}$ are isometric with normals $n_{\pm}$, we
argue as above to see that $W_{\pm}$ are bundles of
isotropic $3$-planes for $\kappa$ which are permuted by $S$.
It follows that exactly one of them is of the form
$\R^4\wedge L$ for some line bundle $L\leq\underline{\R}^4$
and so, after perhaps passing to a double cover of $\Sigma$,
of the form $\R^4\wedge x$ for some map
$x\colon\Sigma\to S^3$.  Without loss of generality, take
$W_+=\R^4\wedge x$ so that, with $F:=F_++F_-$, we have
\begin{equation*}
  \d F=\omega\wedge x,\qquad n_{+}+n_-=\sqrt{2}n\wedge x
\end{equation*}
with $n,\omega$ $\R^4$-valued and orthogonal to $x$.  Since
$n_{\pm}$ are unit normals to $F_{\pm}$, we rapidly conclude
that $n$ has unit length, $n\wedge x$ is normal to $F$ and
$\omega$ is orthogonal to $n$.  Moreover, computing
$(\II_+-\II_-)\vol$ yields
\begin{equation}
  \label{eq:11}
  \omega\wedge x\wedge n\wedge\d x=-2\Re(iq)\vol.
\end{equation}
In particular, since $\Re(iq)$ is non-zero and therefore
non-degenerate, $\d x$ injects so that $x$ immerses.  We
need to show that $x$ is isothermic with retraction form
$\eta=\omega\wedge(x+\q)$, for which it suffices to show
that \eqref{eq:7} holds.  For this, we have
\begin{equation*}
  0=\d^2F=\d\omega\wedge x+\omega\wedge\d x
\end{equation*}
where the first summand takes values in $\R^4\wedge L$ and
the second in $\Wedge^2L^{\perp}$ so that each vanishes
separately.  Thus $\d\omega$ is $L$-valued.  On the other
hand, from the vanishing of $\omega\wedge\d x$, we deduce
that $\im\omega=\im\d x$ so that, in particular, $\omega$ is
orthogonal to $x$ (and $n$ is a normal to $x$ in $S^3$).
Further, from \eqref{eq:10} and \eqref{eq:11}, we get
\begin{equation*}
  \ip{\omega,\d x}=2\Re q
\end{equation*}
so that, in particular, $\ip{\omega,\d x}$ is symmetric.
With this in hand,
\begin{equation*}
  \ip{\d\omega,x}=\d\ip{\omega,x}+\ip{\omega\wedge\d x}=0,
\end{equation*}
since the first summand vanishes by $\ip{\omega,x}=0$ and
the second by symmetry of $\ip{\omega\wedge \d x}$.  We
conclude that $\d\omega=0$ and we are done.

\section{Quaternionic formalism}
\label{sec:quat-form}

Let us compare this analysis with that of
Kamberov--Pedit--Pinkall \cite{KamPedPin98} which employs a
quaternionic formalism.

We therefore view $\R^4$ as the algebra $\H$ of quaternions
and $\R^3\cong\fso(3)$ as the Lie algebra of imaginary
quaternions $\Im\H$ with commutator as Lie bracket.  The
inner product on $\R^4$ is now given by
\begin{equation*}
  \ip{a,b}=\Re(a\bar{b})=\Re(\bar{a}b).
\end{equation*}
The Lie algebra $\fso(4)\cong\Im\H\oplus\Im\H$ with the
latter acting on $\H$ by
\begin{equation}
  \label{eq:12}
  (z_{L},z_{R})c=z_{L}c-cz_{R}.
\end{equation}
To compare \eqref{eq:12} with \eqref{eq:1}, we calculate:
\begin{align*}
  (a\wedge b)c&=\ip{a,c}b-a\ip{b,c}\\
  &=\half\left(a\bar{c}b+c\bar{a}b-a\bar{b}c-a\bar{c}b\right)\\
  &=\half\left(
    c\bar{a}b-a\bar{b}c\right)=\half\left(c\Im(\bar{a}b)-\Im(a\bar{b})c\right),\\
  \intertext{since $\Re(\bar{a}b)=\Re(a\bar{b})$,}
 &=\half(\Im(b\bar{a}),-\Im(\bar{a}b))c.
\end{align*}
We conclude that we have an isomorphism $\Wedge^2\H\cong \Im\H\oplus\Im\H$:
\begin{equation}
  \label{eq:13}
  a\wedge b\mapsto \half\left(\Im(b\bar{a}),-\Im(\bar{a}b)\right)
\end{equation}
We can now express our construction as follows: start with
an isothermic $x\colon\Sigma\to S^3\sub\H$.  Then there is
a closed $\H$-valued $1$-form $\omega$, orthogonal to $x$, with
\begin{equation*}
  \bar{\omega}\wedge\d x=0=\d x \wedge\bar{\omega},
\end{equation*}
where we use quaternionic multiplication to multiply
coefficients in the wedge products.

Now our Bonnet pair $F_{\pm}$ satisfy
\begin{equation*}
  \d F_+=\half x\bar{\omega},\qquad \d F_-=-\half\bar{\omega}x,
\end{equation*}
where we note that the right hand sides are already
imaginary since $\Re(x\bar{\omega})=0$.
\begin{rem}
  In fact, our $F_-$ differs from that in \cite{KamPedPin98}
  by a sign.
\end{rem}


\end{document}